\documentclass[11pt]{amsart}
\usepackage[margin=1.1in]{geometry}
\usepackage{amsmath,amssymb,amsthm}
\usepackage{booktabs}
\usepackage{xcolor}
\usepackage{hyperref}
\hypersetup{colorlinks=true,linkcolor=blue!60!black,citecolor=blue!60!black,urlcolor=blue!60!black}

\newtheorem{theorem}{Theorem}
\newtheorem{lemma}[theorem]{Lemma}

\newtheorem{corollary}[theorem]{Corollary}
\newtheorem{conjecture}[theorem]{Conjecture}
\theoremstyle{definition}
\newtheorem{definition}[theorem]{Definition}
\newtheorem{example}[theorem]{Example}
\theoremstyle{remark}

\newcommand{\seqnum}[1]{\href{https://oeis.org/#1}{\rm \underline{#1}}}

\newcommand{\Fac}{\mathcal F}
\newcommand{\U}{\mathcal U}

\title{The M\"obius function of factorizations with a simple prime factor}

\author{Madjid Mirzavaziri}
\address{Department of Pure Mathematics, Ferdowsi University of Mashhad, Mashhad, Iran}
\email{mirzavaziri@um.ac.ir, mirzavaziri@gmail.com}

\author{Daniel Yaqubi}
\address{Department of Computer Science, University of Torbat-e Jam, Torbat-e Jam, Iran}
\email{yaqubi@tjamcaas.ac.ir, daniel\_yaqubi@yahoo.es}
\thanks{Corresponding author: Daniel Yaqubi.}


\begin{document}
\begin{abstract}
The poset of factorizations of a positive integer $n$ (ordered by refinement) is known to have a very irregular Möbius function. We prove that this function becomes completely explicit whenever $n$ is divisible by a prime $q$ to the first power only ($q^2 \nmid n$). Specifically, between a factorization with $k$ factors and $n$ itself, the Möbius function is given by $(-1)^{k-1}$ multiplied by the number of ways to order the factors that are not divisible by $q$. As a consequence, the sum of the absolute values of the Möbius function over this interval equals twice the number of ordered factorizations of $n/q$.
\end{abstract}

\maketitle

\section{Introduction}
Writing a positive integer $n$ as a product of integers greater than one (ignoring order), we can order these factorizations by refinement: one factorization refines another if it is obtained by splitting some of the factors further. This relation turns the set of factorizations of $n$ into a finite poset, with the prime factorization at the bottom and the trivial factorization $n$ at the top. When $n$ is squarefree, this poset is simply the lattice of set partitions, whose Möbius function was computed independently by Schützenberger, Frucht, and Rota; in particular, the Möbius value for an $r$-element set is $(-1)^{r-1}(r-1)!$ \cite{Sch, Fru, Rot}. However, as soon as any prime factor repeats, the poset ceases to be a lattice, and its Möbius function becomes notoriously difficult to control.

This difficulty arises even for prime powers. The factorizations of $p^k$ correspond to the integer partitions of $k$, and Ziegler showed that this poset fails to be Cohen--Macaulay for large $k$ \cite{Ziegler}, ruling out standard topological tools for computing its Möbius function. For a general integer $n$, the poset is equivalent to the partition poset of the multiset of prime factors of $n$, studied by Hersh and Kleinberg \cite{HK}. They observed that the classical Möbius function is not multiplicative across blocks and instead introduced a multiplicative deformation $\mu'$, proving a product formula for it when one element of the multiset appears with multiplicity one. Meanwhile, the sum of absolute values of the ordinary Möbius function over all factorizations of $n$ was recently added to the On-Line Encyclopedia of Integer Sequences as \seqnum{A369713} \cite{OEIS}, but without an explicit formula.

The purpose of this note is to show that the standard Möbius function itself satisfies a product formula analogous to that of Hersh and Kleinberg and much more whenever $n$ possesses at least one simple prime factor. Our main result is stated as follows.

\begin{theorem}\label{thm:main}
Let $n=mq$ where $q$ is a prime not dividing $m$, and let $x$ be a factorization of $n$ into $k$ factors. Let $c_1,\dots,c_s$ be the multiplicities of the distinct factors of $x$ that are not divisible by $q$. Then
\[
\mu(x,n)=(-1)^{k-1}\,\frac{(k-1)!}{c_1!\cdots c_s!}.
\]
In particular, if $m=p_1^{e_1}\cdots p_r^{e_r}$ and $\Omega(m)=e_1+\dots+e_r$, then the Möbius function between the prime factorization of $n$ and $n$ itself is
\[
(-1)^{\Omega(m)}\,\frac{\Omega(m)!}{e_1!\cdots e_r!}.
\]
\end{theorem}

The number on the right is the number of ways to arrange the prime factors of $m$ in a row, so its absolute value counts the ordered factorizations of $m$ into primes. The theorem cannot be proved by a naive induction on upper intervals, because an upper interval of a factorization poset is in general not again a factorization poset; we explain this in Example~\ref{ex:collision}. The way around it is to work in a slightly larger class of posets, which is closed under taking upper intervals, and to find a sign-reversing argument there.

\section{Coarsening posets}

It is convenient to encode a positive integer by its vector of exponents. A factorization of $n=p_1^{a_1}\cdots p_t^{a_t}$ is then a finite multiset of nonzero vectors in $\mathbb N^t$ whose sum is $(a_1,\dots,a_t)$, and it refines another factorization if its vectors can be grouped so that the group sums are the vectors of the other. We generalize this as follows.

\begin{definition}
Let $W$ be a finite multiset of nonzero vectors in $\mathbb N^t$. A \emph{coarsening} of $W$ is a multiset obtained by partitioning $W$ into groups and replacing each group by its sum. The set $\U(W)$ of coarsenings of $W$ is ordered by declaring $y\le z$ if $z$ is a coarsening of $y$.
\end{definition}

The relation is transitive because groupings compose, and antisymmetric because a proper coarsening has fewer vectors, so $\U(W)$ is a poset with least element $W$ and greatest element the one-vector multiset $\{\Sigma W\}$. If $W$ consists of unit vectors, $e_i$ taken $a_i$ times, then $\U(W)$ is the poset $\Fac(n)$ of factorizations of $n$. The point of the definition is the following observation, immediate from the definitions.

\begin{lemma}\label{lem:upper}
For every $y\in\U(W)$, the upper interval $[y,\{\Sigma W\}]$ of $\U(W)$ is equal to $\U(y)$.
\end{lemma}

\begin{example}\label{ex:collision}
For $n=p^3q^3r^2$ consider the factorization $x=p\cdot p\cdot q\cdot q\cdot r\cdot r\cdot pq$. Its seven factors are of four kinds, with multiplicities $2,2,2,1$, and one might expect the interval $[x,n]$ to look like the factorization poset of an integer with exponents $2,2,2,1$, whose Möbius function between bottom and top is $90$. In fact $\mu(x,n)=92$. The reason is that the factor $pq$ already present in $x$ can also be produced by merging a $p$ with a $q$, so that different groupings of the factors of $x$ give the same coarsening. The interval $[x,n]$ is $\U(x)$, but it is not isomorphic to a factorization poset.
\end{example}

Let us call a vector of $W$ \emph{marked} if its last coordinate is positive, and say that $W$ is \emph{pointed} if exactly one of its vectors is marked and that vector has last coordinate $1$. In a coarsening of a pointed multiset there is again exactly one marked vector, namely the sum of the group containing the marked vector of $W$, so every element of $\U(W)$ is pointed. For a pointed multiset $y$ we write $|y|$ for its number of vectors and $N(y)$ for the number of distinct sequences formed by its unmarked vectors; if these have multiplicities $c_1,\dots,c_s$, then $N(y)=(|y|-1)!/(c_1!\cdots c_s!)$.

\begin{lemma}\label{lem:key}
If $W$ is pointed and $|W|\ge2$, then
\[
\sum_{y\in\U(W)}(-1)^{|y|-1}N(y)=0 .
\]
\end{lemma}

\begin{proof}
Let $w$ be the marked vector of $W$ and $W'=W\setminus\{w\}$, which is nonempty. For $y\in\U(W)$ let $\rho(y)$ be the multiset obtained from $y$ by subtracting $w$ from its marked vector and discarding the result if it is zero. Since the marked vector of $y$ is $w$ plus the sum of a group of vectors of $W'$, the multiset $\rho(y)$ is a coarsening of $W'$. Conversely, let $z\in\U(W')$ have distinct vectors $v_1,\dots,v_h$ with multiplicities $d_1,\dots,d_h$. The elements $y$ with $\rho(y)=z$ are exactly $z\cup\{w\}$ and, for $j=1,\dots,h$, the multiset obtained from $z$ by replacing one copy of $v_j$ with $v_j+w$; these $h+1$ multisets are distinct, since their marked vectors are distinct. For $y=z\cup\{w\}$ we have $|y|=|z|+1$ and $N(y)=|z|!/(d_1!\cdots d_h!)$. For the $j$-th element we have $|y|=|z|$ and
\[
N(y)=\frac{(|z|-1)!}{d_1!\cdots(d_j-1)!\cdots d_h!}=\frac{d_j}{|z|}\cdot\frac{|z|!}{d_1!\cdots d_h!}.
\]
Hence the contribution of the fibre of $z$ is
\[
(-1)^{|z|}\frac{|z|!}{d_1!\cdots d_h!}\Bigl(1-\sum_{j=1}^h\frac{d_j}{|z|}\Bigr)=0,
\]
because $d_1+\dots+d_h=|z|$. Summing over $z\in\U(W')$ gives the claim.
\end{proof}

\begin{theorem}\label{thm:pointed}
If $W$ is pointed, then the Möbius function of $\U(W)$ between its least and greatest elements is $(-1)^{|W|-1}N(W)$.
\end{theorem}

\begin{proof}
We argue by induction on $|W|$; for $|W|=1$ both sides equal $1$. Let $|W|\ge2$ and write $\hat1$ for the greatest element. By the defining recursion of the Möbius function,
\[
\mu(W,\hat1)=-\sum_{\substack{y\in\U(W)\\ y\ne W}}\mu(y,\hat1).
\]
For $y\ne W$ the interval $[y,\hat1]$ equals $\U(y)$ by Lemma~\ref{lem:upper}, and $y$ is pointed with $|y|<|W|$, so $\mu(y,\hat1)=(-1)^{|y|-1}N(y)$ by induction. Lemma~\ref{lem:key} now gives $\mu(W,\hat1)=(-1)^{|W|-1}N(W)$.
\end{proof}

\begin{proof}[Proof of Theorem~\ref{thm:main}]
Take the exponent vectors with respect to the primes of $n$, listing $q$ last. Since $q$ divides $n$ exactly once, a factorization $x$ of $n$ is a pointed multiset, its marked vector being the factor divisible by $q$. By Lemma~\ref{lem:upper} the interval $[x,n]$ of $\Fac(n)$ is $\U(x)$, and Theorem~\ref{thm:pointed} gives $\mu(x,n)=(-1)^{k-1}N(x)$, which is the first formula. For the prime factorization the unmarked vectors are $e_i$ with multiplicity $e_i$, and $k-1=\Omega(m)$.
\end{proof}

\section{The total absolute Möbius function}

Let $a(n)=\sum_{x}|\mu(x,n)|$, the sum running over all factorizations $x$ of $n$; this is the sequence \seqnum{A369713} \cite{OEIS}. Let $H(m)$ be the number of ordered factorizations of $m$, that is, of ways to write $m$ as an ordered product of integers greater than one, with $H(1)=1$. It satisfies $H(m)=\sum_{d\mid m,\,d<m}H(d)$ for $m>1$, and its Dirichlet series is $1/(2-\zeta(s))$ \cite{Hille}.

\begin{corollary}\label{cor:total}
If $n=mq$ with $q$ a prime not dividing $m$ and $m>1$, then $a(n)=2H(m)$.
\end{corollary}

\begin{proof}
By Theorem~\ref{thm:main}, $|\mu(x,n)|$ is the number of orderings of the factors of $x$ not divisible by $q$. A factorization together with such an ordering is the same thing as a divisor $d$ of $m$, recording the factor $qd$, together with an ordered factorization of $m/d$. Hence $a(n)=\sum_{d\mid m}H(m/d)=H(m)+\sum_{d\mid m,\,d<m}H(d)=2H(m)$.
\end{proof}

For instance $a(60)=2H(12)=16$ and $a(24)=2H(8)=8$, in agreement with \cite{OEIS}. The entry \cite{OEIS} also contains the conjecture that $a(n)\le n^c$ for some constant $c$. In the case covered by the corollary this holds with $c=2$: by induction, $H(m)\le\sum_{d\mid m,\,d<m}d^2\le m^2\sum_{j\ge2}j^{-2}<m^2$, so $a(n)\le 2(n/q)^2\le n^2/2$. The Dirichlet series shows that the exponent can be lowered to any number larger than the root $\rho\approx1.7286$ of $\zeta(\rho)=2$.

In the course of these computations, carried out with two independent programs, we obtained $a(16)=a(81)=4$ and $a(36)=a(64)=10$, whereas the entry \cite{OEIS} lists $5$, $5$, $11$ and $11$. In all four cases the poset is small enough to check by hand; for $n=16$ the Möbius values of $16$, $8\cdot2$, $4\cdot4$, $4\cdot2\cdot2$ and $2\cdot2\cdot2\cdot2$ relative to $16$ are $1,-1,-1,1,0$. These $n$ have no simple prime factor, and all of them admit factorizations with two equal composite factors.

\section{Integers without a simple prime factor}

When every prime divides $n$ at least twice, Example~\ref{ex:collision} shows that the structure is genuinely different, and we have no general formula. Table~\ref{tab:values} lists the Möbius function between the prime factorization and $n$ for some exponent patterns, computed exhaustively.

\begin{table}[ht]
\centering
\begin{tabular}{lr@{\qquad}lr@{\qquad}lr}
\toprule
exponents & $\mu$ & exponents & $\mu$ & exponents & $\mu$\\
\midrule
$(2,2)$ & $-1$ & $(3,3)$ & $-3$ & $(4,4)$ & $-8$\\
$(3,2)$ & $2$ & $(4,3)$ & $5$ & $(5,4)$ & $13$\\
$(4,2)$ & $-3$ & $(5,3)$ & $-8$ & $(6,4)$ & $-22$\\
$(5,2)$ & $4$ & $(6,3)$ & $13$ & $(5,5)$ & $-18$\\
$(6,2)$ & $-5$ & $(2,2,2)$ & $-16$ & $(3,2,2)$ & $31$\\
$(7,2)$ & $7$ & $(8,2)$ & $-9$ & $(4,2,2)$ & $-52$\\
$(3,3,2)$ & $-72$ & $(2,2,2,2)$ & $-312$ & & \\
\bottomrule
\end{tabular}
\medskip
\caption{The Möbius function between the prime factorization and $n$ for some exponent patterns without a simple prime.}
\label{tab:values}
\end{table}

For prime powers $p^k$ the value is $1$ and $-1$ for $k=1,2$ and, according to our computations, $0$ for $3\le k\le 10$. The column $n=p^aq^2$ is deceptive: it reads $-1,2,-3,4,-5$ for $2\le a\le6$, suggesting $(-1)^{a+1}(a-1)$, but continues with $7$ and $-9$. The signs, on the other hand, behave regularly.

\begin{conjecture}
The Möbius function between the prime factorization of $n$ and $n$ is either $0$ or has the sign $(-1)^{\Omega(n)+1}$, and it vanishes for $n=p^k$ with $k\ge3$.
\end{conjecture}

For integers with a simple prime factor the sign statement is part of Theorem~\ref{thm:main}. For the others it holds in every case we computed, which includes all $n$ with $\Omega(n)\le8$ and all entries of Table~\ref{tab:values}.

\section*{Acknowledgement}
Computations and parts of the exposition were carried out with the assistance of the AI system Claude (Anthropic).

\end{document}